\documentclass{amsart}

\usepackage{amssymb}
\usepackage{amsfonts}
\usepackage{dsfont}
\usepackage{mathrsfs}
\usepackage{enumitem}

\newcommand{\ncmd}{\newcommand}

\ncmd{\btheo}{\begin{theo}$\!\!\!$ -- }
\ncmd{\etheo}{\end{theo}}
\ncmd{\bmtheo}{\begin{maintheo}$\!\!\!$ -- }
\ncmd{\emtheo}{\end{maintheo}}
\ncmd{\bfait}{\begin{fait}$\!\!\!$ -- }
\ncmd{\efait}{\end{fait}}
\ncmd{\bpro}{\begin{prop}$\!\!\!$ -- }
\ncmd{\epro}{\end{prop}}
\ncmd{\bpreu}{{\sc Proof --}\ }
\ncmd{\epreu}{$\;\;\;\square$}
\ncmd{\bdefi}{\begin{defi}$\!\!\!$ -- }
\ncmd{\edefi}{\end{defi}}
\ncmd{\bco}{\begin{cor}$\!\!\!$ -- }
\ncmd{\eco}{\end{cor}}
\ncmd{\ble}{\begin{lem}$\!\!\!$ -- }
\ncmd{\ele}{\end{lem}}
\ncmd{\bno}{\begin{nota}$\!\!\!$ -- }
\ncmd{\eno}{\end{nota}}
\ncmd{\bre}{\begin{rem}$\!\!\!$ --  \begin{em}}
\ncmd{\ere}{\end{em} \end{rem}}
\ncmd{\bque}{\begin{que}$\!\!\!$ -- \begin{em}}
\ncmd{\eque}{\end{em} \end{que}}
\ncmd{\bconj}{\begin{conj}$\!\!\!$ -- \begin{em}}
\ncmd{\econj}{\end{em} \end{conj}}
\ncmd{\bexe}{\begin{exe}$\!\!\!$ -- \begin{em}}
\ncmd{\eexe}{\end{em} \end{exe}}

\ncmd{\C}{\mathds{C}}
\ncmd{\R}{\mathds{R}}
\ncmd{\F}{\mathds{F}}
\ncmd{\N}{\mathds{N}}
\ncmd{\Q}{\mathds{Q}}
\ncmd{\Z}{\mathds{Z}}
\ncmd{\Am}{\mathscr{A}}
\ncmd{\Bm}{\mathscr{B}}
\ncmd{\Cm}{\mathscr{C}}
\ncmd{\Dm}{\mathscr{D}}
\ncmd{\Em}{\mathscr{E}}
\ncmd{\Fm}{\mathscr{F}}
\ncmd{\Gm}{\mathscr{G}}
\ncmd{\Hm}{\mathscr{H}}
\ncmd{\Ims}{\mathscr{I}}
\ncmd{\Lm}{\mathscr{L}}
\ncmd{\Mm}{\mathscr{M}}
\ncmd{\Nm}{\mathscr{N}}
\ncmd{\Pm}{\mathscr{P}}
\ncmd{\Rm}{\mathscr{R}}
\ncmd{\Sm}{\mathscr{S}}
\ncmd{\Tm}{\mathscr{T}}
\ncmd{\Um}{\mathscr{U}}
\ncmd{\Fc}{\mathfrak{F}}
\ncmd{\Mcc}{\mathfrak{M}}
\ncmd{\Tc}{\mathfrak{T}}
\ncmd{\Ur}{$U_{0,r}$-group}
\ncmd{\Urs}{$U_{0,r}$-subgroup}
\ncmd{\Urp}{$U_{0,r}$-groups}
\ncmd{\Ursp}{$U_{0,r}$-subgroups}
\ncmd{\hUr}{homogeneous $U_{0,r}$-group}
\ncmd{\hUrs}{homogeneous $U_{0,r}$-subgroup}
\ncmd{\hUrp}{homogeneous $U_{0,r}$-groups}
\ncmd{\hUrsp}{homogeneous $U_{0,r}$-subgroups}
\ncmd{\wU}{$\widetilde{U}$}
\ncmd{\wwU}{\widetilde{U}}
\ncmd{\wV}{$\widetilde{V}$}
\ncmd{\wK}{\widetilde{K}}
\ncmd{\Ir}{{\rm Im\hspace{0.05cm}}}
\ncmd{\Kr}{{\rm Ker\hspace{0.05cm}}}
\ncmd{\GL}{{\rm GL\hspace{0.05cm}}}
\ncmd{\PSL}{{\rm PSL}}
\ncmd{\SO}{{\rm SO}}
\ncmd{\Autd}{{\rm Aut}_{def}}
\ncmd{\Auta}{{\rm Aut}_{alg}}
\ncmd{\Aut}{{\rm Aut}}
\ncmd{\Stab}{{\rm Stab}}
\ncmd{\rk}{{\rm rk\hspace{0.05cm}}}
\ncmd{\ad}{{\rm ad}}

\ncmd{\bi}{\begin{itemize}}
\ncmd{\ei}{\end{itemize}}
\ncmd{\be}{\begin{enumerate}}
\ncmd{\ee}{\end{enumerate}}

\ncmd{\dis}{\displaystyle}
\ncmd{\ov}{\overline}
\ncmd{\nn}{\noindent}
\ncmd{\gui}{\textquotedblleft}
\ncmd{\ichap}{\^{\i}}

\newtheorem{lem}{Lemma}[section]
\newtheorem{defi}[lem]{Definition}
\newtheorem{theo}[lem]{Theorem}
\newtheorem{maintheo}[lem]{Main Theorem}
\newtheorem{prop}[lem]{Proposition}
\newtheorem{cor}[lem]{Corollary}
\newtheorem{nota}[lem]{Notation}
\newtheorem{rem}[lem]{Remark}
\newtheorem{fait}[lem]{Fact}
\newtheorem{que}[lem]{Question}
\newtheorem{conj}[lem]{Conjecture}
\newtheorem{exe}[lem]{Example}

\begin{document}

\title{Bad groups in the sense of Cherlin}
\author{Olivier Fr\'econ}
\address{Laboratoire de Math\'ematiques et Applications, 
Universit\'e de Poitiers, 
T\'el\'eport 2 - BP 30179, Boulevard Marie et Pierre Curie, 
86962 Futuroscope Chasseneuil Cedex, France}
\email{olivier.frecon@math.univ-poitiers.fr}
\date{\today}
\begin{abstract}
There exists no bad group (in the sense of Gregory Cherlin), 
namely any simple group of Morley rank 3 is isomorphic to 
${\rm PSL}_2(K)$ for an algebraically closed field $K$.
\end{abstract}
\maketitle


\section{Introduction}

Model theory is a branch of mathematical logic 
concerned with the study of classes of mathematical structures 
by considering first-order sentences and formulas. 
The {\em Morley rank} is a model-theoretical notion of abstract dimension. 
It generalizes the dimension of an algebraic variety 
(when the ground field is algebraically closed). 
There are other notions of abstract dimension, 
the importance of the Morley rank lies on {\em Morley's Categoricity Theorem} below, 
which \gui can be thought of as the beginning of modern model theory''
(David Marker \cite[p. 2]{Marker}) and the following Baldwin and Zilber Theorems. 

We remember that a {\em theory} is a set of first-order $\Lm$-sentences for a language $\Lm$, 
it is {\em complete} if for any sentence $\phi$, 
either $\phi$ or $\neg\phi$ belongs to $T$, 
and a theory is {\em $\kappa$-categorical} for some 
cardinal $\kappa$ if, up to isomorphism, it has exactly one model of cardinality $\kappa$ 
(cf. \cite[Chapters 1 and 2]{Marker} for more details). 

\bfait
Let $T$ be a complete theory in a countable language.
\begin{description}[font= $\bullet$ \normalfont \rm]
\item {\em (Morley's Categoricity Theorem, \cite{Morley})} 
If $T$ is {$\kappa$-categorical} for some 
uncountable $\kappa$, then $T$ is $\kappa$-categorical for every uncountable $\kappa$.
\item {\em (Baldwin, \cite{Bal73})} 
If $T$ is uncountably categorical, then it is {\em of finite Morley rank}.
\item {\em (Zilber, \cite{Zil77})} The theory of an {\em infinite simple group} of finite Morley rank is uncountably categorical. 
\end{description}
\efait

In this paper, we are concerned with {\em groups} of finite Morley rank. 
The main example of such a group is an algebraic group defined over an algebraically closed field 
in the field language (Zilber, \cite{Zil77}). 
In the late seventies, Gregory Cherlin \cite[\S 6]{Che79} and Boris Zilber \cite{Zil77} 
formulated independently the following algebraicity conjecture. 

\bconj
{\em (Cherlin-Zilber Conjecture {\rm or} Algebraicity Conjecture)}
An infinite simple group of finite Morley rank is algebraic over an algebraically closed field.
\econj

This is the main conjecture on groups of finite Morley rank, and it is still open. 
Most of studies on groups of finite Morley rank focus on this conjecture. 
Actually, the original Cherlin Conjecture concerned simple {\em $\omega$-stable} groups, 
but the substantial litterature on the Algebraicity Conjecture treats only the finite Morley rank case.

The Algebraicity Conjecture has been proved for several important classes of groups 
including locally finite groups \cite{Thomas}. 
The main theorem on groups of finite Morley rank ensures that 
any simple group of finite Morley rank with an infinite abelian subgroup of exponent 2 
satisfies the Cherlin-Zilber Conjecture \cite{ABC}.

However, in despite of  numerous papers on the subject, 
the Cherlin-Zilber Conjecture is still open, even for groups of Morley rank 3. 
As a matter of fact, in  \cite{Che79}, the Algebraicity Conjecture 
was formulated as a result from an analysis of simple groups of Morley rank 3. 
The main result of \cite{Che79} can be summarized as follows, 
where a {\em bad group} is a nonsolvable group of Morley rank 3 containing 
no definable subgroup of Morley rank 2.

\bfait\label{thcherlin}
{\em (Cherlin, \cite{Che79})}
Let $G$ be an infinite simple group of Morley rank at most 3. 
Then $G$ has Morley rank 3, and one of the following two assertions is satisfied:
\bi
\item there is an algebraically closed field $K$ such that $G\simeq \PSL_2(K)$,
\item $G$ is a {bad group}.
\ei
\efait

Thus bad groups became a major obstacle to the Cherlin-Zilber Conjecture. 
These groups have been studied in \cite{Che79,Nes89rk3} and \cite{Nes91}, 
whose results are summarized in Facts \ref{bad} and \ref{factNes} respectively. 
Later, it was shown that no bad group is 
existentially closed \cite{JO04} or linear \cite{MusPoi06}. 
However, these groups appeared  very resistant, and 
very sparse other information was known on bad groups. 

Furthermore, Nesin has shown in \cite{Nes91} 
that a bad group acts on a natural geometry, 
which is not very far from being a non-Desarguesian projective plane of Morley rank 2. 
However, Baldwin discovered non-Desarguesian projective planes of Morley rank 2 \cite{Bal94}. 
Thus, the question of the existence, or not, of a bad group was still fully open. 
In this paper, we show that bad groups do not exist. 

\bmtheo
There is no bad group.
\emtheo

Note other more general notions of bad groups 
have been introduced independently by Corredor \cite{Cor89} and 
by Borovik and Poizat \cite{BP90}, where a {\em bad group} is defined to be a nonsolvable connected group 
of finite Morley rank all of whose proper connected definable subgroups are nilpotent. 
Such a bad group has similar properties 
to original bad groups. 
Moreover, later Jaligot will introduce a more general notion of  bad groups \cite{Jal01}, 
and he will obtain similar results. 
However, we recall that, in this paper, a  {\em bad group} is defined to be 
nonsolvable, of Morley rank 3, and containing no definable subgroup of Morley rank 2.

\medskip

Our proof of Main Theorem goes as follows. 
First we note that it is sufficient to study {\em simple} bad groups 
since for any bad group $G$, the quotient group $G/Z(G)$ is a {\em simple} bad group 
by \cite[\S 4, Introduction]{Nes89rk3}. 

Then we fix a simple bad group $G$, and 
we introduce a notion of lines 
as cosets of Borel subgroups of $G$ (Definition \ref{defiline}). 
In \S \ref{secline}, we study their behavior, 
mainly in regards with conjugacy classes of elements of $G$. 

In \S \ref{secplane}, we propose a definition of a plane 
(Definition \ref{defiplane}). 
This section is dedicated to prove that $G$ contains a plane (Theorem \ref{thplane}). 
This result is the key point of our demonstration. 
Roughly speaking, we show that for each nontrivial element $g$ of $G$ 
such that $g=[u,v]$ for $(u,v)\in G\times G$, 
the union of the preimages of $g$, by maps of the form ${\rm ad}_v:G\to G$ 
defined by ${\rm ad}_v(x)=[x,v]$, is almost a plane, and from this, we obtain a plane.

In last section \S \ref{secfin}, 
we try to show that our notions of lines and planes 
provide a structure of projective space over the group $G$. 
Indeed, such a structure would provide a division ring 
(see \cite[p. 124, Theorem 7.15]{Harts}), 
and probably it would be easy to conclude. 
However, a contradiction occurs along the way, and achieves our proof.

{\em \underline{Note :} in a very recent preprint \cite{PoiWag}, 
by analyzing the present paper, Poizat and Wagner 
generalize our main result to other groups, 
and they eliminate other groups of Morley rank.}

\subsection*{The other simple groups of dimension 3}

\bi
\item If $G$ is a non-bad simple group of Morley rank 3, 
then $G$ is isomorphic to $\PSL_2(K)$ 
for an algebraically closed field $K$ (Fact \ref{thcherlin}). 
As in \S \ref{secline}, we may define a {\em line} in $G$ 
to be a coset of a connected subgroup of dimension 1, 
and we may define a plane as in \S \ref{secplane}. 
It is possible to show that two sorts of planes occur: 
the cosets of Borel subgroups, 
and the subsets of the form $aJ$ where $J$ is defined to be 
\bi
\item the set of involutions when the characteristic $c$ of $K$ is not 2;
\item $J=\{j\in G~|~j^2=1\}$ when $c=2$.
\ei
%
The plane $J$ is normalized by $G$, 
and there is no such a plane in a bad group (Lemma \ref{mapalphaaut}). 
Another important difference between $G$ and a bad group 
is to be the presence of a Weyl group. 
Indeed, the first lemma of this paper is not verified in $G$ (Lemma \ref{lem1}), 
because we have $jT=Tj$ for any torus $T$ and any involution $j\in N_G(T)\setminus T$. 

\item The group $\SO_3(\R)$ is not of finite Morley rank, or even stable \cite{Nes89rk3}. 
However, our definitions of lines and  planes naturally extend to $\SO_3(\R)$. 
Then, as above, the set $J$ of involutions in $\SO_3(\R)$ forms a plane, 
and the presence of a Weyl group is again a major difference between 
$\SO_3(\R)$ and bad groups. Moreover, we note that the plane $J$ 
has a structure of projective plane, whereas this is false in $\PSL_2(K)$ 
\cite[Fact 8.15]{bn1}.
\ei

\section{Background material}

A thorough analysis of groups of finite Morley rank can be found in \cite{bn1} and \cite{ABC}. 
In this section we recall 
some definitions and known results. 

\subsection{Borovik-Poizat axioms}

Let $(G,\,\cdot\,,^{-1},1,\cdots)$ be a group equipped with additional structure. 
This group $G$ is said to be {\em ranked} if there is a function 
\gui {\rm rk}'' which assigns to each nonempty definable set $S$ an integer, 
its \gui dimension'' $\rk(S)$, and which satisfies the following axioms 
for every definable sets $A$ and $B$.
{\em\begin{description}[font=\normalfont \rm]
\item[Definition] For any integer $n$, 
$\rk(A)>n$ if and only if $A$ contains an infinite family 
of disjoint definable subsets $A_i$ of rank $n$.
\item[{Definability}] For any uniformly definable family $\{A_b~:~b\in B\}$ 
of definable sets, and for any $n\in \N$, the set $\{b\in B~:~\rk(A_b)=n\}$ 
is also definable.
\item[{Finite Bounds}] For any uniformly definable family $\Fc$ of finite subsets of $A$, 
the sizes of the sets in $\Fc$ are bounded.
\end{description}}
It is shown in \cite{poigrsta} that the groups $(G,\,\cdot\,,\cdots)$ as above 
satisfy a fourth axiom, namely the {\em additivity axiom}, 
and they are precisely the groups of finite Morley rank. 
Moreover, the function $\rk$ 
assigns to each definable set its Morley rank. 
In this paper, as in \cite{bn1} and \cite{ABC}, the Morley rank 
will be denoted by $\rk$. 

\subsection{Morley degree}

A nonempty definable set $A$ is said to have {\em Morley degree 1} 
if for any definable subset $B$ of $A$, either $\rk B<\rk A$ or $\rk(A\setminus B)<\rk A$. 
The set $A$ is said to have {\em Morley degree} $d$ if $A$ is the disjoint union 
of $d$ definable sets of Morley degree 1 and Morley rank $\rk A$. 

\bfait
\bi
\item\cite[Lemmas 4.12 and 4.14]{bn1} Every nonempty definable set has a unique degree.
\item \cite[Proposition 4.2]{bn1} Let $X$ and $Y$ be definable subsets of Morley degree $d$ and $d'$ respectively. 
Then $X\times Y$ has Morley degree $dd'$. 
\item \cite[\S 2.2]{Che79} A group of finite Morley rank has Morley degree 1 if and only if it is {\em connected}, 
namely it has no proper definable subgroup of finite index.
\ei
\efait

Moreover, the following elementary result will be useful for us.

\bfait\label{deg1}
Let $f:E\to F$ be a definable map. If the set $E$ has Morley degree 1 
and $r=\rk f^{-1}(y)$ is constant for $y\in F$, 
then the Morley degree of $F$ is 1. 
\efait

\bpreu
Let $B$ be a definable subset of $F$ of Morley rank $\rk F$. 
We show that $\rk(F\setminus B)<\rk F$. 
By the additivity axiom, we have $\rk E=r+\rk F$ 
and $$\rk f^{-1}(B)=r+\rk B=r+\rk F=\rk E$$
Since $E$ has Morley degree 1, we obtain $\rk f^{-1}(F\setminus B)=\rk(E\setminus f^{-1}(B))<\rk E$, 
and by the additivity axiom again, 
$$\rk(F\setminus B)=\rk f^{-1}(F\setminus B)-r<\rk E-r=\rk F$$
so $F$ has Morley degree 1.
\epreu

\subsection{Bad groups}

Main properties of bad groups are summarized in the following facts, 
where a {\em Borel subgroup} of a bad group $G$  
is defined to be an infinite definable proper subgroup of $G$.

\bfait\label{bad}
{\em (\cite[\S 5.2]{Che79} and \cite{Nes89rk3})}
Let $G$ be a simple bad group, and $B$ be a Borel subgroup of $G$. 
\be
\item $B=C_G(b)$ for $b\in G\setminus \{1\}$,
\item $B$ is connected, abelian, self-normalizing and of Morley rank 1,
\item $C_G(x)$ is a Borel subgroup for each nontrivial element $x$ of $G$,
\item if $A$ is another Borel subgroup of $G$, then $A$ is conjugate with $B$,
and either $A=B$ or $A\cap B=\{1\}$,
\item $G=\bigcup_{g\in G}B^g$, 
\item $G$ has no involution.
\ee
\efait

\bfait\label{factNes}
{\em \cite[Lemma 18]{Nes91}}
Let $A$ and $B$ be two distinct Borel subgroups of a simple bad group $G$. 
Then $\rk(ABA)=3$, $\rk(AB)=2$, and $AB$ has Morley degree 1. 
\efait

The following result is due to Delahan and Nesin, 
it is proved for a more general notion of bad groups, 
and it is used in our final argument.

\bfait\label{invoautobad}
{\em \cite[Proposition 13.4]{bn1}}
A simple bad group $G$ cannot have an involutive definable automorphism.
\efait

\section{Lines}\label{secline}

In this paper, $G$ denotes a fixed {simple bad group}. 
We fix a Borel subgroup $B$ of $G$ and we denote by $\Bm$ the set of Borel subgroups of $G$. 

In this section, we define a {\em line} of $G$, and we provide their basic properties. 
We note that, by conjugation of Borel subgroups (Fact \ref{bad} (4)), any Borel subgroup 
is a {\em line} in the following sense.

\bdefi\label{defiline}
A {\em line} of $G$ is a subset of the form $uBv$ for two elements $u$ and $v$ of $G$.

We denote by $\Lambda$ the set of lines of $G$. 
\edefi

We note that, by Fact \ref{bad} (2), each line has Morley rank 1 and Morley degree 1.



\ble\label{lem1}
Let $uBv$ and $rBs$ be two lines. 
Then $uBv=rBs$ if and only if $uB=rB$ and $Bv=Bs$.
\ele

\bpreu
We may assume that $uBv=rBs$. Then we have 
$$B=u^{-1}rBsv^{-1}=u^{-1}rsv^{-1}B^{sv^{-1}}$$
so $u^{-1}rsv^{-1}\in B^{sv^{-1}}$ and $B=B^{sv^{-1}}$. 
Now $sv^{-1}$ belongs to $B$ since $B$ is self-normalizing by Fact \ref{bad}. 
Hence we obtain $Bv=Bs$, and the equality $uB=rB$ follows from $uBv=rBs$. 
\epreu

\medskip

By the above lemma, the set $\Lambda$ identifies with $(G/B)_l\times (G/B)_r$ 
where $(G/B)_l$ (resp. $(G/B)_r$) denotes the set of left cosets (resp. right cosets) 
of $B$ in $G$. 
Then $\Lambda$ is a definable set. 
Moreover, since $G$ is connected of Morley rank 3 
and $B$ has Morley rank 1, the Morley 
rank of $\Lambda$ is 4 and its Morley degree is 1. 
In particular, $\Lambda$ is a uniformly definable family. 



\ble\label{lem3}
Two distinct elements $x$ and $y$ of $G$ lie in one and only one line $l(x,y)$. 
Moreover, the map $l:\{(x,y)\in G\times G~|~x\neq y\}\to \Lambda$ is definable.
\ele

\bpreu
By Fact \ref{bad} (5), there exists $v\in G$ such that $y^{-1}x$ belongs to $B^v$. 
Then $x$ and $y$ lie in $uBv$ for $u=yv^{-1}$.

Now, if $rBs$ is a line containing $x$ and $y$, then we find two elements $b_1$ and $b_2$ of $B$ 
such that $x=rb_1s$ and $y=rb_2s$. 
Thus $y^{-1}x=s^{-1}b_2^{-1}b_1s$ is a nontrivial element of $B^s$. 
But $y^{-1}x$ belongs to $B^v$ by the choice of $v$, 
hence we have $B^s=B^v$ (Fact \ref{bad} (4)). 
Since $B$ is self-normalizing, $sv^{-1}$ belongs to $B$ and we obtain $Bs=Bv$, 
so there exists $b\in B$ such that $s=bv$. 
This implies that $u=yv^{-1}=(rb_2s)(s^{-1}b)=rb_2b$ belongs to $rB$, 
and $rBs=uBv$ is the unique line containing $x$ and $y$.

Moreover, since $\Lambda$ is a uniformly definable family, 
the set $\{(x,y)\in G\times G~|~x\neq y\}\times \Lambda$ is definable, and 
$$\Gamma=\{((x,y),uBv)\in (G\times G)\times \Lambda~|~
x\neq y,~ x\in uBv,~y\in uBv\}$$
is a definable subset of it. 
But $\Gamma$ is precisely the graph of the map $l$, 
hence $l$ is definable.
\epreu

\ble\label{lem2}
If $uBv=(uBv)^g$ for $uBv\in \Lambda\setminus\Bm$ and $g\in G$, 
then $g=1$.
\ele
 
\bpreu
We have $uBv=g^{-1}uBvg$, so $uB=g^{-1}uB$ and $Bv=Bvg$ by Lemma \ref{lem1}, 
and $g$ belongs to the Borel subgroups $B^{u^{-1}}$ and $B^v$. 
If $g$ is nontrivial, then  $B^{u^{-1}}=B^v$ (Fact \ref{bad} (4)), 
and $vu$ belongs to $N_G(B)=B$. 
Consequently $u$ belongs to $v^{-1}B$, and we obtain $uBv=B^v$, contradicting $uBv\not\in \Bm$. 
Thus $g=1$.
\epreu



\bdefi
For each $g\in G$ and each definable subset $X$ of $G$, 
we consider the following subsets of $\Lambda$:
$$\Lm(g,X)=\{l(g,x)\in\Lambda~|~x\in X\setminus\{g\}\}$$
$$\Lambda_X=\{\lambda\in \Lambda~|~\lambda\cap X\ {\rm is\ infinite}\}$$
\edefi

Since the map $l$ is definable (Lemma \ref{lem3}), the set $\Lm(g,X)$ 
is definable for each $g\in G$ and each definable subset $X$ of $G$. 
Moreover, by Definablity axiom, the set 
$\Lambda_X=\{\lambda\in \Lambda~|~\rk(\lambda\cap X)=1\}$ is definable too.

\ble\label{lemdegfin}
Let $\lambda_1,\ldots,\lambda_n$ be $n$ lines. 
Then $\lambda_1\cup\cdots\cup \lambda_n$ is a definable set of Morley rank 1 
and Morley degree $n$. 
\ele

\bpreu
For each $i$, the set $A_i=\lambda_i\cap(\bigcup_{j\neq i}\lambda_j)$ has at most $n-1$ elements 
by Lemma \ref{lem3}, and $\lambda_1\cup\cdots\cup \lambda_n$ is the disjoint union 
of $\lambda_1\setminus A_1,\ldots,\lambda_n\setminus A_n,\cup_{i=1}^nA_i$. 
Since each line $\lambda_i$ has Morley rank 1 
and Morley degree 1 (Fact \ref{bad} (2)), the result follows.
\epreu

\ble\label{lemdef12}
If $\Lambda_0$ is a definable subset of $\Lambda$, 
then $\bigcup\Lambda_0$ is a definable subset of $G$. 
Moreover, if $\Lambda_0$ is infinite, 
then $\bigcup\Lambda_0$ has Morley rank at least 2.
\ele

\bpreu
Since $\Lambda_0$ is a definable subset of 
the uniformly definable family $\Lambda$, 
the set $\bigcup\Lambda_0=\{x\in G~|~\exists \lambda\in \Lambda_0,~x\in \lambda\}$ 
is definable. 
Moreover, if $\Lambda_0$ is infinite, 
then $\bigcup\Lambda_0$ has Morley rank at least 2 by Lemma \ref{lemdegfin}.
\epreu

\bco\label{cordef12}
The subset $\bigcup\Lambda_X$ of $G$ is definable for each definable subset $X$ of $G$. 
\eco

\section{Planes}\label{secplane}

Our aim is to find a definable structure of projective space 
on our bad group $G$. 
In this section, we introduce a notion of planes, 
and we show that $G$ has such a plane (Theorem \ref{thplane}). 
We fix a definable subset $X$ of $G$, of Morley rank 2.

\bdefi\label{defiplane}
The definable subset $X$ of $G$ is said to be a {\em plane} if it satisfies $\ov{X}=X$ where 
$$\ov{X}=\{g\in G~|~{\rk}(\Lm(g,X))=1\}$$
\edefi

\ble\label{procomm1}
The set $\ov{X}$ is a definable subset of $\bigcup\Lambda_X$. 
\ele

\bpreu
If $g\in G$ does not belong to $\bigcup\Lambda_X$, 
then $l(g,x)\cap X$ is finite for each $x\in X$, 
and since $X$ has Morley rank 2, 
the set $\Lm(g,X)$ has Morley rank 2, so $g\not\in\ov{X}$. 
Thus $\ov{X}$ is contained in $\bigcup\Lambda_X$. 

We show that $\ov{X}$ is definable. 
We consider the set $$A=\{(g,\lambda)\in G\times \Lambda~|~g\in \lambda,~\exists x\in X\setminus\{g\},~x\in \lambda\}$$ 
and the map $f:A\to G$ defined by $f(g,\lambda)=g$. 
We note that, since $\Lambda$ is a uniformly definable family, 
$A$ is definable, 
and $f$ is definable too. 
Moreover, the preimage by $f$ of each $g\in G$ 
is $f^{-1}(g)=\{g\}\times \Lm(g,X)$, 
and we have $\rk (f^{-1}(g))=\rk\Lm(g,X)$. 
Consequently, we obtain $\ov{X}=\{g\in G~|~\rk (f^{-1}(g))=1\}$, 
and $\ov{X}$ is definable.
\epreu

\ble\label{lemrkdm2}
The Morley ranks of $\Lambda_X$ and of $\bigcup\Lambda_X$ are at most 2. 
Moreover, $\Lambda_X$ is infinite if and only if $\rk\bigcup\Lambda_X=2$.
\ele

\bpreu
We consider the surjective definable map 
$$l_0:(X\times X)\cap l^{-1}(\Lambda_X)\to \Lambda_X$$ defined by $l_0(x,y)=l(x,y)$. 
For each $\lambda\in \Lambda_X$, we have 
$l_0^{-1}(\lambda)=\{(x,y)\in (\lambda\cap X)\times (\lambda\cap X)~|~x\neq y\}$, 
and since $\rk \lambda=1$, we obtain $\rk(\lambda\cap X)=1$ 
and $\rk l_0^{-1}(\lambda)=2$. 
But we have 
$$\rk ((X\times X)\cap l^{-1}(\Lambda_X))\leq \rk(X\times X)=2\rk X=4$$
hence $\rk \Lambda_X$ is at most $4-2=2$. 

We show that $\rk\bigcup\Lambda_X\leq 2$. 
We consider the definable set 
$$A=\{(x,\lambda)\in G\times \Lambda_X~|~x\in \lambda\setminus X\}$$
and the definable map $l_1:A\to \Lambda_X$ defined by $l_1(x,\lambda)=\lambda$. 
For each $\lambda\in \Lambda_X$, we have $\rk \lambda=1=\rk(\lambda\cap X)$, 
so, since each line has Morley degree 1, the preimage $l_1^{-1}(\lambda)$ is finite. Consequently 
we obtain $\rk A\leq \rk\Lambda_X\leq 2$. 
But the definable map $l_2:A\to (\bigcup\Lambda_X)\setminus X$, defined by $l_2(x,\lambda)=x$, 
is surjective, hence the Morley rank of $(\bigcup\Lambda_X)\setminus X$ is at most $\rk A\leq 2$. 
Since $X$ has Morley rank 2, we obtain $\rk\bigcup\Lambda_X\leq 2$.

Now it follows from Lemmas \ref{lemdegfin} and \ref{lemdef12} 
that $\Lambda_X$ is infinite if and only if $\rk\bigcup\Lambda_X=2$.
\epreu

\bpro\label{lemdm1}
For each $g\in \ov{X}$, we have ${\rk}(\Lm(g)\cap \Lambda_X)=1$.

Moreover, if $X$ has Morley degree 1, then $\ov{X}=\{g\in G~|~{\rk}(\Lm(g)\cap \Lambda_X)=1\}$ 
and $G\setminus \ov{X}=\{g\in G~|~\Lm(g)\cap \Lambda_X~{\rm is\ finite}\}$.
\epro

\bpreu
First we note that $\Lm(g)\cap \Lambda_X=\Lm(g,X)\cap \Lambda_X$ for any $g\in G$. 
For each $g\in G$, we consider the definable map 
$l_g:X\setminus\{g\}\to \Lm(g,X)$ defined by $l_g(x)=l(g,x)$. 
In particular, the preimage $l_g^{-1}(\lambda)$ of each 
$\lambda\in \Lm(g,X)$ is $(\lambda\cap X)\setminus\{g\}$. 

We show that ${\rk}(\Lm(g,X)\cap \Lambda_X)\leq 1$ for each $g\in G$. 
We may assume that $\Lambda_X$ is infinite. 
Then, by Lemma \ref{lemrkdm2}, the set $\cup\Lambda_X$ has Morley rank 2. 
Let $g\in G$ and $u_g:\cup(\Lm(g)\cap\Lambda_X)\setminus\{g\}\to \Lm(g)\cap\Lambda_X$ 
be the map defined by $u_g(x)=l(g,x)$. 
Since each line has Morley rank 1, 
the preimage of each element of $\Lm(g)\cap\Lambda_X$ has Morley rank 1. 
Consequently, we have $$\rk (\Lm(g)\cap\Lambda_X)=
\rk\cup(\Lm(g)\cap\Lambda_X)-1\leq \rk\cup\Lambda_X-1=1$$

Let $g\in \ov{X}$. We show that ${\rk}(\Lm(g)\cap \Lambda_X)=1$. 
For each $\lambda\in \Lm(g,X)\setminus \Lambda_X$, 
the set $l_g^{-1}(\lambda)=(\lambda\cap X)\setminus\{g\}$ is finite, and 
since $g\in\ov{X}$, we have $\rk\Lm(g,X)=1$. 
Consequently,  $l_g^{-1}( \Lm(g,X)\setminus \Lambda_X)$ has Morley rank at most 1, 
and $l_g^{-1}( \Lm(g,X)\cap\Lambda_X)$ has Morley rank $\rk X=2$. 
But the set $l_g^{-1}(\lambda)=(\lambda\cap X)\setminus\{g\}$ is infinite of Morley rank 1 
for each $\lambda\in \Lm(g,X)\cap \Lambda_X$. 
Hence we obtain $\rk (\Lm(g,X)\cap \Lambda_X)=2-1=1$.

Now we assume that $X$ has Morley degree 1. Let $g\in G$ 
such that ${\rk}(\Lm(g,X)\cap \Lambda_X)=1$. We show that $g\in\ov{X}$. 
Since the set $l_g^{-1}(\lambda)=(\lambda\cap X)\setminus\{g\}$ is infinite of Morley rank 1 
for each $\lambda\in \Lm(g,X)\cap \Lambda_X$, 
the set $l_g^{-1}( \Lm(g,X)\cap\Lambda_X)$ has Morley rank 
$$1+{\rk}(\Lm(g,X)\cap \Lambda_X)=2=\rk X$$ 
Then, since $X$ has Morley degree 1, 
the preimage of $\Lm(g,X)\setminus \Lambda_X$ has Morley rank at most 1. 
Moreover, for each $\lambda\in \Lm(g,X)\setminus \Lambda_X$, 
the preimage $l_g^{-1}(\lambda)=(\lambda\cap X)\setminus\{g\}$ is finite and non-empty, 
so we obtain $$\rk(\Lm(g,X)\setminus \Lambda_X)=\rk l_g^{-1}(\Lm(g,X)\setminus \Lambda_X)\leq 1$$
This shows that $\rk\Lm(g,X)=1$ and $g\in\ov{X}$. 

Furthermore, since ${\rk}(\Lm(g,X)\cap \Lambda_X)\leq 1$ for each $g\in G$, 
we obtain $G\setminus \ov{X}=\{g\in G~|~\Lm(g)\cap \Lambda_X~{\rm is\ finite}\}$, 
as desired.
\epreu

\bco\label{corY2}
We have $\rk(\ov{X}\setminus X)\leq 1$.
\eco

\bpreu
We remember that $\ov{X}$ is definable by Lemma \ref{procomm1}, 
so the sets $Y=\ov{X}\setminus X$ and $A=\{(y,\lambda)\in Y\times \Lambda_X~|~y\in \lambda\}$ 
are definable too. 
Let $l_Y:A\to Y$ and $l_D:A\to \Lambda_X$ be the definable maps defined by 
$l_Y(y,\lambda)=y$ and $l_D(y,\lambda)=\lambda$ respectively. 
On the one hand, for each $\lambda\in\Lambda_X$, 
the set $\lambda\cap X$ is infinite, and since $\lambda$ has Morley rank 1 
and Morley degree 1 (Fact \ref{bad} (2)), the set $\lambda\cap Y$ is finite 
and $l_D^{-1}(\lambda)$ has Morley rank at most 0. 
This implies $\rk A\leq\rk\Lambda_X\leq 2$ (Lemma \ref{lemrkdm2}). 
On the other hand, for each $y\in Y$, we have ${\rk}(\Lm(y)\cap \Lambda_X)=1$ 
by Proposition \ref{lemdm1}, so $l_Y^{-1}(y)$ has Morley rank 1, 
and we obtain $\rk A=1+\rk Y$. 
Consequently, the Morley rank of $Y$ is at most 1.
\epreu

\ble\label{lemmaxdl}
For each $g\in G$, the set $\Lm(g,X)$ is infinite. 
\ele

\bpreu
Indeed, $\cup\Lm(g,X)$ is definable (Lemma \ref{lemdef12}) 
and contains $X$.  Since $\rk X=2$, we obtain $\rk(\cup\Lm(g,X))\geq 2$, 
and $\Lm(g,X)$ is infinite (Lemma  \ref{lemdegfin}).
\epreu

\bco\label{planedeg1}
If the Morley degree of $X$ is not 1, then $\rk\ov{X}<2$. 
In particular, any plane has Morley degree 1.
\eco

\bpreu
Let $n$ be the Morley degree of $X$, and $X_1,\ldots,X_n$ 
be $n$ definable subsets of $X$ of Morley rank 2 and Morley degree 1 
such that $X$ is the disjoint union of $X_1,\ldots,X_n$. 
For each $g\in\ov{X}$, we have $\rk\Lm(g,X)=1$, 
so we obtain $\rk\Lm(g,X_i)\leq 1$ for each $i$, 
and $g\in\ov{X_i}$ for each $i$ by Lemma \ref{lemmaxdl}. 
Thus $\ov{X}$ is contained in $\ov{X_1}\cap\ov{X_2}$. 
Since $X_1\cap X_2=\emptyset$, the set $\ov{X}$ is contained in 
$(X_1\cap Y_2)\cup (Y_1\cap X_2)\cup (Y_1\cap Y_2)$ where 
$Y_1=\ov{X_1}\setminus X_1$ and $Y_2=\ov{X_2}\setminus X_2$. 
Since $Y_1$ and $Y_2$ have Morley rank at most 1 by Corollary \ref{corY2}, 
we obtain $\rk\ov{X}<2$.
\epreu

\ble\label{plangen1}
We assume that $X$ has Morley degree 1, 
and that $Y$ is another definable subset of $G$ of Morley rank 2 and Morley degree 1. 
If $X\cap Y$ has Morley rank 2, then $\ov{X}=\ov{Y}$.
\ele

\bpreu
Let $g\in G$. 
If $g$ belongs to $\ov{X\cap Y}$, then we have $\rk\Lm(g,X\cap Y)=1$. 
Since $X$ has Morley degree 1 and $X\cap Y$ 
has Morley rank 2, the set $X\setminus Y$ has Morley rank at most 1, 
and the set $\Lm(g,X\setminus Y)$ has Morley rank at most 1. 
Thus $\Lm(g,X)$ has Morley rank 1, and $g$ belongs to $\ov{X}$. 

Conversely, if $g\in\ov{X}$, then $\Lm(g,X)$ has Morley rank 1, 
so $\Lm(g,X\cap Y)\subseteq \Lm(g,X)$ has Morley rank at most 1. 
Then Lemma \ref{lemmaxdl} gives $g\in\ov{X\cap Y}$. 
This shows that $\ov{X\cap Y}=\ov{X}$. 
By the same way, we obtain $\ov{X\cap Y}=\ov{Y}$, so $\ov{X}=\ov{Y}$.
\epreu


\medskip

For each $a\in G$, let $\Lm(a)=\Lm(a,G)$ be the (definable) set of lines containing $a$. 
Moreover, we note that $\Lm(1)=\Bm$.

\ble\label{lemrk21}
Let $\Lambda_0$ be a definable subset of $\Lambda$. 
If $\rk\cup\Lambda_0=2$, 
then we have $\rk(\Lm(g)\cap\Lambda_0)\leq 1$ for each $g\in G$. 

Moreover, if further $\rk\Lambda_0=2$, then 
the set $\{g\in G~|~\rk(\Lm(g)\cap\Lambda_0)=1\}$ has Morley rank 2.
\ele

\bpreu
We show that $\rk(\Lm(g)\cap\Lambda_0)\leq 1$ for each $g\in G$. 
Let $g\in G$ and $l_g:\cup(\Lm(g)\cap\Lambda_0)\setminus\{g\}\to \Lm(g)\cap\Lambda_0$ 
be the map defined by $l_g(x)=l(g,x)$. 
Since each line has Morley rank 1, 
the preimage of each element of $\Lm(g)\cap\Lambda_0$ has Morley rank 1. 
Consequently, we have $$\rk (\Lm(g)\cap\Lambda_0)=
\rk\cup(\Lm(g)\cap\Lambda_0)-1\leq \rk\cup\Lambda_0-1=1$$
as desired. 

We suppose further that $\rk\Lambda_0=2$, and 
we show that $\{g\in G~|~\rk(\Lm(g)\cap\Lambda_0)=1\}$ has Morley rank 2. 
Let $U=\cup\Lambda_0$, $A=\{(u,\lambda)\in U\times \Lambda_0~|~u\in \lambda\}$ and 
$f:A\to \Lambda_0$ be the map defined by $f(u,\lambda)=\lambda$. 
Then $A$ and $f$ are definable, and the preimage $f^{-1}(\lambda)$ 
of each $\lambda\in \Lambda_0$ has Morley rank $\rk \lambda=1$, 
so $\rk A=1+\rk\Lambda_0=3$. 
Now let $h:A\to U$ be the map defined by $h(u,\lambda)=u$. 
It is a definable map, and the preimage $h^{-1}(u)$ 
of each $u\in U$ has Morley rank either 0, or 1 by the previous paragraph. 

But the preimage of $U_0=\{u\in U~|~\rk h^{-1}(u)=0\}$ has Morley rank 
$$\rk h^{-1}(U_0)=\rk U_0\leq \rk U=2<\rk A$$
so the preimage of $U_1=\{u\in U~|~\rk h^{-1}(u)=1\}$ has Morley rank 3. 
Hence we obtain $\rk U_1=3-1=2$. 
Moreover, we note that 
$$U_1=\{u\in U~|~\rk(\Lm(u)\cap\Lambda_0)=1\}=\{g\in G~|~\rk(\Lm(g)\cap\Lambda_0)=1\}$$
so $\{g\in G~|~\rk(\Lm(g)\cap\Lambda_0)=1\}$ has Morley rank 2.
\epreu

\bpro\label{prodeg1}
Let $X$ be a definable subset of $G$ of Morley rank 2 and Morley degree 1. 
Then $\rk\ov{X}=2$ if and only if $\Lambda_X$ has Morley rank 2. 

In this case, $\Lambda_X$ and $\ov{X}$ have Morley degree 1, 
and $\ov{X}$ contains a generic definable subset of $X$.
\epro

\bpreu
We consider the definable set 
$A=\{(x,\lambda)\in \ov{X}\times \Lambda_X~|~x\in \lambda\}$  
and the definable maps $l_1:A\to \ov{X}$ and $l_2:A\to \Lambda_X$ defined 
by $l_1(x,\lambda)=x$ and $l_2(x,\lambda)=\lambda$ respectively. 
By Proposition \ref{lemdm1}, the preimage $l_1^{-1}(g)$ of each element $g$ of $\ov{X}$ 
has Morley rank 1, so $\rk A=1+\rk\ov{X}$. 
Moreover, the preimage $l_2^{-1}(\lambda)$ of each $\lambda\in\Lambda_X$ 
has Morley rank at most 1, so $\rk A\leq 1+\rk \Lambda_X$. 
Then we obtain $\rk\ov{X}\leq \rk\Lambda_X$. 
In particular, it follows from Lemma \ref{lemrkdm2} that if $\rk\ov{X}=2$, 
then $\rk\Lambda_X=2$. 
Hence we may assume that $\rk\Lambda_X=2$. 

At this stage, Lemma \ref{lemrkdm2} gives $\rk\cup\Lambda_X=2$, 
and by Lemma \ref{lemrk21} and Proposition \ref{lemdm1}, we obtain $\rk\ov{X}=2$. 
Moreover, it follows from Corollary \ref{corY2} that $\ov{X}$ has Morley degree 1 
and that $X\cap\ov{X}$ is a generic definable subset of $X$ contained in $\ov{X}$.

We show that the Morley degree of $\Lambda_X$ is 1.
Let $l_0:\{(x,y)\in X\times X~|~x\neq y\}\to \Lambda$ 
be the definable map defined by $l_0(x,y)=l(x,y)$. 
Since the Morley degree of $X$ is 1, 
the one of $\{(x,y)\in X\times X~|~x\neq y\}$ is 1 too.
For each 
$\lambda\in\Lambda_X$, we have $\rk l_0^{-1}(\lambda)=\rk((\lambda\cap X)\times (\lambda\cap X))=2$. 
Since $\rk\Lambda_X=2$, we obtain 
$$\rk l_0^{-1}(\Lambda_X)=2+\rk\Lambda_X=4=\rk\{(x,y)\in X\times X~|~x\neq y\}$$
and since the Morley degree of $\{(x,y)\in X\times X~|~x\neq y\}$ is 1, 
the Morley degree of $l_0^{-1}(\Lambda_X)$ is 1 too. 
Now the Morley degree of $\Lambda_X$ is 1 by Fact \ref{deg1}.
\epreu

\ble\label{linederive}
Let $g$ be a nontrivial element such that $g=[u,v]$ 
for $(u,v)\in G\times G$. 
Then we have $\{x\in G~|~[x,v]=g\}=C_G(v)u$ and $\{y\in G~|~[u,y]=g\}=C_G(u)v$. 
In particular, they are two lines and have Morley rank 1 and Morley degree 1.
\ele

\bpreu
The equalities are obvious. Moreover, 
by Fact \ref{bad}, the sets $C_G(v)u$ and $C_G(u)v$ are two lines, 
and they have Morley rank 1 and Morley degree 1. 
\epreu

\ble\label{lem4}
For each $a\in G$, the set $a^G\cap B$ has exactly one element.
\ele

\bpreu
We may assume $a\neq 1$. 
By Fact \ref{bad} (5), there is $g\in G$ such that $a^g$ belongs to $B$. 
If $a^h\in B$ for $h\in G$, then $a$ is a nontrivial element of $B^{g^{-1}}\cap B^{h^{-1}}$. 
By Fact \ref{bad} (4), we obtain $B^{g^{-1}}=B^{h^{-1}}$, 
and $h^{-1}g$ belongs to $N_G(B)=B$. 
But $B$ is abelian (Fact \ref{bad} (2)), so $h^{-1}g$ centralizes $a^h$, 
and $a^h=(a^h)^{h^{-1}g}=a^g$. 
Hence $a^G\cap B=\{a^g\}$.
\epreu

\medskip

The following result isolates a step of the proof of Theorem \ref{thplane}. 
Its proof and the one of Theorem \ref{thplane} were originally a lot more complicated, 
and Bruno Poizat provided a simplification.

For each $g\in G$, we consider the following definable subset of $G$:
$$X(g)=\{x\in G~|~\exists y\in G,~[x,y]=g\}$$

\bpro\label{proX2}
For each nontrivial element $g$ of $G$, 
the set $X(g)$ has Morley rank at most 2. 
\epro

\bpreu
We assume toward a contradiction that $X(g)$ has Morley rank 3. 
Then the Morley rank of $X(g^z)$ is 3 for each $z\in G$. 
We recall that, by Fact \ref{bad}, the conjugacy class $g^G$ of $g$ 
has Morley rank $\rk g^G=\rk G-\rk C_G(g)=2$

We consider $V=\{(x,y)\in G\times G~|~[x,y]\in g^G\}$ and 
the definable surjective map $f:V\to g^G$ defined by $f(x,y)=[x,y]$. 
For each $z\in G$, we have $$f^{-1}(g^z)=\{(x,y)\in G\times G~|~[x,y]=g^z\}$$
and by Lemma \ref{linederive}, this set has Morley rank 
$\rk f^{-1}(g^z)=\rk X(g^z) +1=3 +1=4$, so $\rk V=4+\rk g^G=6$. 
Since $G\times G$ is a connected group of Morley rank 6, 
the set $V$ is a definable generic subset of $G\times G$, 
and there is $(x,y)\in V$ such that $(y,x)$ belongs to $V$. 
Thus $[x,y]\in g^G$ and its inverse $[y,x]\in g^G$ are conjugate, 
and they are equal by Lemma \ref{lem4}, contradicting that $G$ has no involution 
(Fact \ref{bad} (6)). 
\epreu

\btheo\label{thplane}
There is a plane in $G$.
\etheo

\bpreu
It is sufficient to show that there is a definable subset $X$ of $G$ 
satisfying the following properties:
\be
\item its Morley rank is 2 and its Morley degree is 1,
\item $\Lambda_X$ has Morley rank 2.
\ee
Indeed, by Proposition \ref{prodeg1}, for such a subset $X$, 
the set $\ov{X}$ has Morley rank 2 and Morley degree 1, and 
it contains a generic definable subset of $X$. 
At this stage, Lemma \ref{plangen1} shows that $Y=\ov{X}$ is a plane.

We fix a nontrivial element $g$ such that $g=[u,v]$ for $(u,v)\in G\times G$.

\medskip

{\em 1. For each $x\in X(g)$, there are infinitely many lines 
containing $x$ and contained in $X(g)$.}

\medskip

Since $x$ belongs to $X(g)$, there is $y\in G$ such that $[x,y]=g$. 
We note that, since $g$ is nontrivial, 
$x$ and $y$ are nontrivial and we have $C_G(x)\neq C_G(y)$. 
In particular, $C_G(x)y$ is a line, and it does not contain 1. 
Thus, for each $c\in C_G(x)$, the set $l_c=C_G(cy)x$ is a line, 
and by Lemma \ref{linederive}, we have $[r,cy]=[x,cy]=[x,y]=g$ for each $r\in l_c$. 
So $l_c$ is a line containing $x$ and contained in $X(g)$. 

If $l_c=l_d$ for two elements $c$ and $d$ of $C_G(x)$, 
then we have $C_G(cy)=C_G(dy)$, and $C_G(cy)$ is a line containing $cy$ and $dy$. 
But $C_G(x)y$ is another line containing $cy$ and $dy$, 
and we have $C_G(x)y\neq C_G(cy)$ because $C_G(x)y$ does not contain $1$. 
Hence Lemma \ref{lem3} gives $c=d$, 
and $\{l_c\in\Lambda~|~c\in C_G(x)\}$ is an infinite family of lines 
containing $x$ and contained in $X(g)$.

\medskip

{\em 2. $\rk X(g)=2$.}

\medskip

By 1., the set $X(g)$ contains infinitely many lines, 
so it has Morley rank at least 2 (Lemma \ref{lemdegfin}), 
and by Proposition \ref{proX2}, it has Morley rank 2. 

\medskip

{\em 3. $\rk \Lambda_{X(g)}=2$.}

\medskip

By Lemma \ref{lemrkdm2}, the set $\Lambda_{X(g)}$ has Morley rank at most 2. 
Since $X(g)$ is infinite by 2., for each positive integer $n$ 
we can find $n$ distinct elements $x_1,\ldots,x_n$ in $X(g)$. 
By 1., the set $\Lambda_i=\{\lambda\in \Lambda_{X(g)}~|~x_i\in\lambda\}$ 
is infinite for each $i$. 
We may assume that its Morley rank is 1 for each $i$. 
Then, since there are finitely many lines containing 
two distinct elements among $x_1,\ldots,x_n$ (Lemma \ref{lem3}), 
the union $\cup_{i=1}^n\Lambda_i$ 
has Morley rank 1 and Morley degree at least $n$. 
This implies that $\Lambda_{X(g)}$ does not have Morley rank 1, 
so $\rk \Lambda_{X(g)}=2$.

\medskip

{\em 4. Conclusion.}

\medskip

By 2., the set $X(g)$ has Morley rank 2. 
Let $d$ be its Morley degree. 
Then $X(g)$ is the disjoint union of definable subsets $X_1,\ldots,X_d$ 
of Morley rank 2 and Morley degree 1. 

For each element $\lambda$ of $\Lambda_{X(g)}$, since $\lambda\cap X(g)$ is infinite 
and since $\lambda$ has Morley rank 1 and Morley degree 1, 
there is a unique $i\in\{1,\ldots,d\}$ such that 
$\lambda\cap X_i$ is infinite, that is $\lambda\in \Lambda_{X_i}$. 
Thus, each $\lambda\in\Lambda_{X(g)}$ belongs to a unique 
definable set $\Lambda_{X_i}$ for $i\in\{1,\ldots,d\}$. 
Hence $\Lambda_{X(g)}$ is the disjoint union of 
$\Lambda_{X_1},\ldots,\Lambda_{X_d}$, and there exists $i\in\{1,\ldots,d\}$ 
such that $\rk\Lambda_{X_i}=2$. 
Now the set ${X_i}$ satisfies the conditions (1) and (2) 
of the beginning of our proof, so $\ov{X_i}$ is a plane.
\epreu

\section{A projective space ?}\label{secfin}

In this section, we analyze planes. 
We remember that, by Theorem \ref{thplane}, the group $G$ has a plane, 
and that by Corollary \ref{planedeg1}, any plane has Morley degree 1. 
The initial goal of this section was to show that, if $X$ and $Y$ are two distinct planes, 
then $\Lambda_X\cap \Lambda_Y$ has a unique element. 
However, along the way, we obtain our final contradiction.

\bdefi
For each line $\lambda$, we consider the following subset of $\Lambda$:
$$\Lm(\lambda)=\{m\in \Lambda~|~\lambda\cap m\ {is\ not\ empty}\}$$
\edefi

\ble\label{lemlml}
For any line $\lambda$, 
the set $\Lm(\lambda)$ is definable, it has Morley rank 3 
and Morley degree 1.
\ele

\bpreu
We consider the definable map $f:\lambda\times (G\setminus \lambda)\to \Lm(\lambda)\setminus \{\lambda\}$ 
defined by $f(x,g)=l(x,g)$. 
By Lemma \ref{lem3}, for each $m\in \Lm(\lambda)\setminus \{\lambda\}$, there is a unique 
element $x$ in $\lambda\cap m$. Moreover, for any $g\in G\setminus \lambda$, we have $f(x,g)=m$ 
if and only if $g\in m\setminus\{x\}$. 
Consequently we have $\rk f^{-1}(m)=\rk m=1$, and 
$$\rk\Lm(\lambda)=\rk(\lambda\times (G\setminus \lambda))-1=3$$ 
Furthermore, since $\lambda$ and $G$ have Morley degree 1, 
the Morley degree of $\lambda\times G$ and $\lambda\times (G\setminus \lambda)$ is 1, 
and the Morley degree of $\Lm(\lambda)\setminus \{\lambda\}$ and $\Lm(\lambda)$ is 1 too 
(Fact \ref{deg1}).
\epreu

\ble\label{lemplan1}
Let $X$ be a plane, and $\lambda\in\Lambda_X$. 
Then $\Lm(\lambda)\cap\Lambda_X$ has Morley rank 2.
\ele

\bpreu
Since $\lambda$ belongs to $\Lambda_X$, the set $\lambda\cap X$ is infinite, 
and since $\lambda$ is a line, we have $\rk(\lambda\cap X)=1$. 
We consider the definable set $$\Am=\{(x,m)\in (\lambda\cap X)\times \Lambda_X~|~m\neq \lambda,~x\in m\}$$
and the definable maps $p:\Am\to \lambda\cap X$ and $q:\Am\to\Lambda_X$ 
defined by $p(x,m)=x$ and $q(x,m)=m$ respectively. 
By Proposition \ref{lemdm1}, the set $p^{-1}(x)$ has Morley rank 1 for each $x\in \lambda\cap X$, 
so $\rk\Am=1+\rk(\lambda\cap X)=2$. 

Moreover, each $m\in \Lambda_X\setminus\{\lambda\}$ contains at most one element of $\lambda$ 
(Lemma \ref{lem3}), so $q$ is an injective map and its image has Morley rank $\rk\Am=2$. 
But the image of $q$ is contained in $(\Lm(\lambda)\cap\Lambda_X)\setminus\{\lambda\}$, 
and we have $\rk\Lambda_X\leq 2$ (Lemma \ref{lemrkdm2}), 
hence $\Lm(\lambda)\cap\Lambda_X$ has Morley rank 2.
\epreu

\ble\label{leminterline}
Let $\lambda_1$ and $\lambda_2$ be two distinct lines. 
Then $\Lm(\lambda_1)\cap\Lm(\lambda_2)$ has Morley rank 2 and Morley degree 1.
\ele

\bpreu
Let $A=\{(x,y)\in \lambda_1\times \lambda_2~|~x\not\in\lambda_1,\, y\not\in\lambda_2\}$, 
and let $f:A\to (\Lm(\lambda_1)\cap\Lm(\lambda_2))\setminus\{\lambda_1,\lambda_2\}$ be the 
map defined by $f(x,y)=l(x,y)$. This map is definable and bijective by Lemma \ref{lem3}. 
Since $\lambda_1$ and $\lambda_2$ are two lines, 
the sets $\lambda_1\times \lambda_2$ and $A$ have Morley rank 2 and Morley degree 1, 
and since $f$ is a definable bijection, 
$\Lm(\lambda_1)\cap\Lm(\lambda_2)$ has Morley rank 2 and Morley degree 1.
\epreu

\bpro\label{proDmX}
If $X$ and $Y$ are two distinct planes, 
then $\Lambda_X\cap \Lambda_Y$ has at most one element.
\epro

\bpreu
Suppose toward a contradiction that $\lambda_1$ and $\lambda_2$ 
are two distinct elements of $\Lambda_X\cap \Lambda_Y$. 
By Lemma \ref{lemplan1}, the sets $\Lm(\lambda_1)\cap\Lambda_X$ and $\Lm(\lambda_2)\cap\Lambda_X$ 
have Morley rank 2. But $\Lambda_X$ has Morley rank 2 and Morley degree 1
by Proposition \ref{prodeg1}, 
hence $\Lm(\lambda_1)\cap\Lm(\lambda_2)\cap\Lambda_X$ has Morley rank 2. 
By the same way, $\Lm(\lambda_1)\cap\Lm(\lambda_2)\cap\Lambda_Y$ has Morley rank 2. 
Thus, since $\Lm(\lambda_1)\cap\Lm(\lambda_2)$ has Morley rank 2 and Morley degree 1 
(Lemma \ref{leminterline}), the set $\Lambda_X\cap \Lambda_Y$ has Morley rank 2. 

Since $\Lambda_X\cap \Lambda_Y$ is infinite, 
the set $U=\cup(\Lambda_X\cap \Lambda_Y)$ has Morley rank at least 2 by Lemma \ref{lemdef12}, 
and since $U$ is contained in $\cup\Lambda_X$, its Morley rank is exactly 2 (Lemma \ref{lemrkdm2}). 
Now the set $Z=\{g\in G~|~\rk(\Lm(g)\cap\Lambda_X\cap \Lambda_Y)=1\}$ 
has Morley rank 2 by Lemma \ref{lemrk21}.
But Proposition \ref{lemdm1} says that $Z$ is contained in $X\cap Y$, 
hence $X\cap Y$ has Morley rank 2 
and Lemma \ref{plangen1} gives $X=Y$, a contradiction.
\epreu

\medskip

From now on, we try to show that the set $\Lambda_X\cap \Lambda_Y$ has exactly one element. 
However, the final contradiction will appear earlier.

\bco\label{coraXb}
Let $X$ be a plane and $(a,b)\in G\times G$. 
Then the following assertions are equivalent:
\bi
\item $aXb=X$
\item $a\Lambda_X b=\Lambda_X$
\item $aXb\cap X$ has Morley rank 2.
\ei
\eco

\bpreu
We note that $aXb$ is a plane, and that $a\Lambda_X b=\Lambda_{aXb}$. 
If $aXb\cap X$ has Morley rank 2, 
then $aXb=X$ by Lemma \ref{plangen1}, and 
if $aXb=X$, then we have $a\Lambda_X b=\Lambda_{aXb}=\Lambda_X$. 
Moreover, if $a\Lambda_X b=\Lambda_X$, then we have $\Lambda_{aXb}=\Lambda_X$ 
and $aXb=X$ by Proposition \ref{proDmX}, 
so $aXb\cap X=X$ has Morley rank 2.
\epreu

\medskip

By Fact \ref{factNes}, if $A$ is a Borel subgroup distinct from $B$, 
then $\rk(ABA)=3$. 
The following result is slightly  more general, and its proof is different.

We recall that, if a group $H$ of finite Morley rank acts definably 
on a set $E$, then the {\em stabilizer} of any definable subset $F$ of $E$ is defined to be 
$$\Stab F=\{h\in H~|~\rk((h\cdot F)\Delta F)<\rk (F)\}$$
where $\Delta$ stands for the symmetric difference. 
It is a definable subgroup of $H$ by \cite[Lemma 5.11]{bn1}.

\ble\label{lemrkABC}
Let $A$ and $C$ be two Borel subgroups distinct from $B$. 
Then $\rk(ABC)=3$.
\ele

\bpreu
We consider the action of $G$ on itself by left multiplication. 
Then we have $b\cdot BC=BC$ for each $b\in B$, 
so $B$ is contained in $\Stab(BC)$. 

We assume toward a contradiction that $C$ is contained in $\Stab(BC)$. 
Since $BC$ has Morley rank 2 and Morley degree 1 (Fact \ref{factNes}), 
we have $\rk(cBC\setminus BC)\leq 1$ for each $c\in C$, 
and since $\rk C=1$, we obtain $\rk(CBC\setminus BC)\leq 2$ 
and $\rk(CBC)=2$, contradicting Fact \ref{factNes}. 
Consequently, $C$ is not contained in $\Stab(BC)$, 
and since $\Stab(BC)$ contains $B$, Fact \ref{bad} implies that $\Stab(BC)=B$. 

We assume toward a contradiction that $\rk(ABC)\neq 3$. 
Since $\rk(BC)=2$, we have $\rk(ABC)=2$ and 
$ABC$ is a disjoint union of finitely many definable subsets 
$E_1,\ldots,E_k$ of Morley rank 2 and Morley degree 1. 
For each $a\in A$, the set $aBC$ has Morley rank $\rk (BC)=2$ and Morley degree 1, 
so there exists a unique $i\in\{1,\ldots,k\}$ such that $\rk(aBC\cap E_i)=2$. 
Since $A$ is infinite, there are $i\in\{1,\ldots,k\}$ and 
two distinct elements $a$ and $a'$ of $A$ such that 
$\rk(aBC\cap E_i)=\rk(a'BC\cap E_i)=2$. 
Since $E_i$ has Morley degree 1, the Morley rank of $aBC\cap a'BC$ is 2, 
and we obtain $\rk(a'^{-1}aBC\cap BC)=2$. 
But $BC$ has Morley degree 1, hence $a'^{-1}a$ belongs to $\Stab(BC)=B$. 
Thus $a'^{-1}a$ belongs to $A\cap B=\{1\}$ (Fact \ref{bad} (4)), 
contradicting that $a$ and $a'$ are distinct. 
So  we have $\rk(ABC)=3$, as desired.
\epreu

\bco\label{corBABC}
Let $A$ and $C$ be two distinct Borel subgroups. 
Then $\rk(BA\cap BC)=1$.
\eco

\bpreu
We may assume $A\neq B$ and $C\neq B$. 
By Fact \ref{factNes}, we have 
$$1=\rk B\leq \rk(BA\cap BC)\leq \rk(BA)=2$$
We assume toward a contradiction that $\rk(BA\cap BC)=2$. 
Since $BC$ has Morley rank 2 and Morley degree 1 (Fact \ref{factNes}), 
the set $E=BC\setminus BA$ has Morley rank at most 1. 
Consequently, $EA$ has Morley rank at most $\rk E+\rk A=2$, 
and since $(BA\cap BC)A\subseteq BA$ has Morley rank 2, 
we obtain $\rk(BCA)=\rk(EA\cup(BA\cap BC)A)=2$, 
contradicting that $BCA$ has Morley rank 3 (Lemma \ref{lemrkABC}).
\epreu

\ble\label{lemborstab}
For any plane $X$, we have $BX\neq X$ and $XB\neq X$.
\ele

\bpreu
We assume toward a contradiction that $BX=X$ for a plane $X$. 
Let $x\in X$. Since $X$ is a plane, Proposition \ref{lemdm1} gives 
$\rk(\Lm(x,X)\cap \Lambda_X)=1$, 
so $\Lm(x,X)\cap \Lambda_X$ is infinite. 
But each line containing $x$ has the form $B^ux$ for $u\in G$, 
hence there exist $u\not\in B$ and $v\not\in B$ 
such that $B^u\neq B^v$, and 
such that $B^ux$ and $B^vx$ belong to $\Lm(x,X)\cap \Lambda_X$.
In particular, there is a co-finite subset $S$ of $B$ 
such that $S^ux$ and $S^vx$ are contained in $X$. 

Now, since $BX=X$, the sets $BS^ux$ and $BS^vx$ are contained in $X$. 
By Fact \ref{factNes}, the set $BB^{u}$, and so $Bu^{-1}B$, has Morley rank 2, 
and since $Bu^{-1}(B\setminus S)$ is a finite union of lines, 
the set $Bu^{-1}(B\setminus S)$ has Morley rank 1 (Lemma \ref{lemdegfin}), 
and $Bu^{-1}S$ has Morley rank 2. 
Thus, the sets $BS^ux=Bu^{-1}Sux$ and $BS^vx=Bv^{-1}Svx$ 
are subsets of $X$ of Morley rank 2, 
and since the Morley degree of $X$ is 1, 
the set $BS^ux\cap BS^vx$ has Morley rank 2. 
This implies that $\rk(BB^u\cap BB^v)=2$, contradicting Corollary \ref{corBABC}. 
Now we have $BX\neq X$ and by the same way, we show that $XB\neq X$. 
\epreu

\bco\label{cormultleft}
For any plane $X$, the stabilizer of $X$ for the action of $G$ 
on itself by left multiplication is finite. 
\eco

\bpreu
By Corollary \ref{coraXb}, we have $\Stab X=\{a\in G~|~aX=X\}$. 
If $\Stab X$ is infinite, then it contains a Borel subgroup, contradicting Lemma \ref{lemborstab}.
\epreu

\bpro\label{prodefaXb}
Let $X$ be a plane. 
Then for each plane $Y$, there exist a unique $a\in G$ 
and a unique $b\in G$ such that $Y=aX=Xb$.
\epro

\bpreu
We fix $\alpha\in G$, and we consider the following definable subset of $\Lambda$:
$$A=\{(\lambda_1,\lambda_2)\in \Lambda\times\Lambda~|~\alpha\in\lambda_1\cap\lambda_2,\,\lambda_1\neq\lambda_2\}$$
We show that $A$ has Morley rank 4 and Morley degree 1.  
Let $U=\{(x,y)\in G\times G~|~y\not\in l(x,\alpha)\}$. 
Then $U$ is a generic definable subset of $G\times G$, 
and it has Morley rank 6 and Morley degree 1. 
Let $f:U\to A$ be the definable surjective map defined by 
$f(x,y)=(l(x,\alpha),l(y,\alpha))$. 
Since each line has Morley rank 1, the preimage of each 
$(\lambda_1,\lambda_2)\in A$ has Morley rank $\rk\lambda_1+\rk\lambda_2=2$, 
and the set $A$ has Morley rank $\rk U-2=4$ and Morley degree 1 (Fact \ref{deg1}). 

For each plane $P$, we consider the following definable set
$$A_P=\{(\lambda_1,\lambda_2)\in \Lambda\times\Lambda~|~
\alpha\in\lambda_1\cap\lambda_2,\,\lambda_1\neq\lambda_2,\,
\exists a\in G,\,a^{-1}\lambda_1\in \Lambda_P,\,a^{-1}\lambda_2\in \Lambda_P \}$$
We show that the set $A_X$
is a generic definable subset of $A$. 
Indeed, for each $a\in \alpha X^{-1}$, 
we have $\alpha\in aX$ and $\rk(\Lm(\alpha)\cap\Lambda_{aX})=1$ 
by Proposition \ref{lemdm1}, 
so the definable set 
$$L_{aX}=\{(\lambda_1,\lambda_2)\in \Lambda_{aX}\times\Lambda_{aX}~|~
\alpha\in\lambda_1\cap\lambda_2,\,\lambda_1\neq\lambda_2\}$$
has Morley rank $2\rk(\Lm(\alpha)\cap\Lambda_{aX})=2$. 
But $\alpha X^{-1}$ has Morley rank $\rk X=2$ and 
it follows from Proposition \ref{proDmX} that $L_{aX}\cap L_{bX}=\emptyset$ 
for any two elements $a$ and $b$ of $\alpha X^{-1}$ such that $aX\neq bX$. 
Moreover, for each $a\in \alpha X^{-1}$, there are finitely many elements 
$b\in \alpha X^{-1}$ such that 
$aX=bX$ (Corollary \ref{cormultleft}). 
Hence the set $A_X=\cup_{a\in \alpha X^{-1}}L_{aX}$ has Morley rank $\rk\alpha X^{-1}+2=4$, 
and it is a generic definable subset of $A$. 

By the same way, $A_Y$ is a generic definable subset of $A$, 
so there exists $(\lambda_1,\lambda_2)\in A_X\cap A_Y$. 
Thus there exist two elements $u$ and $v$ of $G$ 
such that two distinct lines $\lambda_1$ and $\lambda_2$ 
belong to $\Lambda_{uX}\cap\Lambda_{vY}$, 
and we obtain $uX=vY$ by Proposition \ref{proDmX}, 
so $Y=aX$ for $a=v^{-1}u$.
By the same way, there exists $b\in G$ such that $Y=Xb$. 

We show the uniqueness of $a$ and $b$. 
Let $S=\{g\in G~|~gX=X\}$. It is a finite subgroup of $G$ 
by Corollary \ref{cormultleft}. 
For each $\alpha\in G$, the previous paragraph gives $\beta\in G$ such that 
$\alpha X=X\beta$. Then, for each $s\in S$, we have $s(\alpha X)=s(X\beta)=X\beta=\alpha X$, 
and we obtain $s^\alpha X=X$ and $s^\alpha\in S$. 
Thus any element $\alpha\in G$ normalizes the finite subgroup $S$, 
and since $G$ is a simple group, $S$ is trivial. 
This proves the uniqueness of $a$, and by the same way we obtain the uniqueness of $b$.
\epreu

\medskip

By the previous result, the set of planes is $\Pm=\{aX~|~a\in G\}$, 
and it identifies with $G$. 
Thus, the set of planes is uniformly definable 
and has Morley rank 3. 



\ble\label{mapalphaaut}
There exists $a\in G$ such that $X^a\neq X$.
\ele

\bpreu
We assume toward a contradiction that $X^a=X$ for each $a\in G$. 
Then for each $uBv\in \Lambda_X$ and each $a\in G$, we have 
$$(uBv)^a\in \Lambda_X^a=\Lambda_{X^a}=\Lambda_X$$
Since $\rk\Lambda_X=2$ (Proposition \ref{prodeg1}), the line $uBv$ is a Borel subgroup 
(Lemma \ref{lem2}), and by conjugacy of Borel subgroups, 
we obtain $\Lambda_X=\Bm$. 
Now we have $\cup\Lambda_X=G$, so $\rk\cup\Lambda_X=3$, 
contradicting Lemma \ref{lemrkdm2}.
\epreu

\medskip

From now on, we are ready for the final contradiction. 
Initially, it was more complicated, 
but Poizat proposed a simplification by introducing the inverted plane. 

\medskip

\bpreu
First we note that for each plane $Y$, the set $y^{-1}Y$ is a plane containing 1, 
and the set $Y^{-1}$ is a plane too. 
We fix a plane $X$ containing 1. 
By Proposition \ref{prodefaXb}, 
there is a bijective map $\mu:G\to G$ defined by $xX=X\mu(x)$, 
and $\mu$ is definable since the set $\Pm$ of planes is uniformly definable. 
Moreover, for each $(a,b)\in G\times G$, 
we have $X\mu(ab)=abX=aX\mu(b)=X\mu(a)\mu(b)$, 
so $\mu$ is an automorphism of $G$. 

Since $X=\ov{X}$ contains $1$, there are infinitely many Borel subgroups in $\Lambda_X$. 
Let $B_1$ and $B_2$ be two distinct Borel subgroups belonging to $\Lambda_X$. 
Then $B_1$ and $B_2$ belong to $\Lambda_{X^{-1}}$ too, 
and we have $X=X^{-1}$ by Proposition \ref{proDmX}. 
By the same way, since the plane $x^{-1}X$ contains 1 for each $x\in X$, 
we have $x^{-1}X=(x^{-1}X)^{-1}=X^{-1}x=Xx$ for each $x\in X$. 
Thus $\mu(x^{-1})=x$ for each $x\in X$, and since $X=X^{-1}$, 
we obtain $\mu^2(x)=x$ for each $x\in X$. 

But $X$ is a definable subset of $G$ of Morley rank 2, 
hence $G$ is generated by $X$, and $\mu$ is an involutive automorphism of $G$. 
Thus $\mu$ is the identity map by Fact \ref{invoautobad}, 
contradicting Lemma \ref{mapalphaaut}.
\epreu

\bre
After Lemma \ref{mapalphaaut} we were ready for a new step to provide 
a structure of projective space over $G$, 
which was the initial goal of our section. 
Indeed, in the first version of this paper, we have shown that, 
if $X$ and $Y$ are two distinct plane, 
then $\Lambda_X\cap \Lambda_Y$ has a unique element. 
\ere

\section*{Acknowledgements}

J'adresse un grand merci \`a Bruno Poizat pour sa lecture attentive de cet article, 
ses commentaires, et plusieurs simplifications de la preuve.

I thank Gregory Cherlin very much for his comments and his interest in this paper, 
Frank O. Wagner for his remarks, and Tuna Alt\i nel for his good advice.

\bibliographystyle{amsplain}
\bibliography{split}

\end{document}